\documentclass[12pt]{article}
\usepackage{enumerate}
\usepackage{latexsym}
\usepackage{amsfonts}
\usepackage[only,ninrm,elvrm,twlrm,sixrm,egtrm,tenrm]{rawfonts}
\usepackage{indentfirst}
\usepackage{amsmath}
\usepackage{graphicx,psfrag}
\usepackage{graphics}
\usepackage{makeidx}
\newtheorem{thm}{Theorem}[section]

\newtheorem{lem}[thm]{Lemma}

\newtheorem{prop}[thm]{Proposition}

\def\qed{\hfill \nopagebreak\rule{5pt}{8pt}}
\def\pf{\noindent {\it Proof.} }
\title{\bf Complete solution for the rainbow numbers of matchings
\footnote {Supported by NSFC and the ``973" project.}}

\author{
\small He Chen, Xueliang Li and Jianhua Tu\\
\small Center for Combinatorics and LPMC, Nankai University\\
\small Tianjin 300071, P.R. China. Email: lxl@nankai.edu.cn\\
}
\date{}

\begin{document}

\maketitle
\begin{abstract}
For a given graph $H$ and $n\geq 1$, let $f(n,H)$ denote the maximum
number $c$ for which there is a way to color the edges of the
complete graph $K_n$ with $c$ colors such that every subgraph $H$ of
$K_n$ has at least two edges of the same color. Equivalently, any
edge-coloring of $K_n$ with at least $rb(n,H)=f(n,H)+1$ colors
contains a rainbow copy of $H$, where a rainbow subgraph of an
edge-colored graph is such that no two edges of it have the same
color. The number $rb(n,H)$ is called the {\it rainbow number of
$H$}. Erd\H{o}s, Simonovits and S\'{o}s showed that $rb(n,K_3)=n$.
In 2004, Schiermeyer used some counting technique and determined the
rainbow numbers $rb(n,kK_2)$ for $k\geq 2$ and $n\geq 3k+3$. It is
easy to see that $n$ must be at least $2k$. So, for $2k \leq
n<3k+3$, the rainbow numbers remain not determined. In this paper we
will use the Gallai-Edmonds structure theorem for matchings to
determine the exact values for rainbow numbers $rb(n,kK_2)$ for all
$k\geq 2$ and $n\geq 2k$, giving a complete solution for the rainbow
numbers of matchings.\\
[0.1in]{\bf Keywords:} edge coloring, rainbow subgraph, rainbow number.\\
[0.1in] {\bf AMS subject classification 2000:} 05C35, 05C55, 05C70.
\end{abstract}

\section{Introduction}

In this paper we consider undirected, finite and simple graphs only,
and use standard notations in graph theory (see \cite{bon} and
\cite{lov}). If $K_n$ is edge-colored and a subgraph $H$ of $K_n$
contains no two edges of the same color, then $H$ is called a {\it
totally multicolored (TMC)} or {\it rainbow subgraph} of $K_n$ and
we say that $K_n$ contains a TMC or rainbow $H$. Let $f(n,H)$ denote
the maximum number of colors in an edge-coloring of $K_n$ with no
TMC $H$. We now define $rb(n,H)$ as the minimum number of colors
such that any edge-coloring of $K_n$ with at least
$rb(n,H)=f(n,H)+1$ colors contains a TMC or rainbow subgraph
isomorphic to $H$. The number $rb(n,H)$ is called the {\it rainbow
number of $H$}.

$f(n,H)$ is called the anti-Ramsey number of $H$. Anti-Ramsey
numbers were introduced by Erd\H{o}s, Simonovits and S\'{o}s in the
1970s. They showed that these are closely related to Tur\'{a}n
numbers. Anti-Ramsey numbers have been studied in
\cite{alon,tao1,erd,axe,mon,tao2,sos} and elsewhere. There are very
few graphs whose anti-Ramsey numbers have been determined exactly.
To the best of our knowledge, $f(n,H)$ is known exactly for large
$n$ only when $H$ is a complete graph, a path, a star, a cycle or a
broom whose maximum degree exceeds its diameter (a broom is obtained
by identifying an end of a path with a vertex of a star) (see
\cite{tao1,ing,mon,tao2,sos}).

For a given graph $H$, let $ext(n,H)$ denote the maximum number of
edges that a graph $G$ of order $n$ can have with no subgraph
isomorphic to $H$. For $H=kK_2$, the values $ext(n,kK_2)$ have been
determined by Erd\H{o}s and Gallai \cite{gal}, where $H=kK_2$ is a
matching $M$ of size $k$.

\begin{thm}(Erd\H{o}s and Gallai \cite{gal}) \
$ext(n,kK_2)=\max\{{2k-1 \choose 2}, {k-1 \choose 2}+(k-1)(n-k+1)\}$
for all $n\geq2k$ and $k\geq1$, that is, for any given graph $G$ of
order $n$, if $|E(G)|>\max\{{2k-1 \choose 2}, {k-1 \choose
2}+(k-1)(n-k+1)\}$, then $G$ contains a $kK_2$, or a matching of
size $k$.
\end{thm}

In 2004, Schiermeyer \cite{ing} used some counting technique and
determined the rainbow numbers $rb(K_n,kK_2)$ for all $k\geq 2$ and
$n\geq 3k+3$.

\begin{thm}(Schiermeyer \cite{ing}) \
$rb(n,kK_2)=ext(n,(k-1)K_2)+2$ for all $k\geq2$ and $n\geq3k+3$.
\end{thm}

It is easy to see that $n$ must be at least $2k$. So, for $2k \leq
n<3k+3$, the rainbow numbers remain not determined. In this paper,
we will use a technique deferent from Schiermeyer \cite{ing} to
determine the exact values of $rb(n,kK_2)$ for all $k\geq 2$ and
$n\geq 2k$. Our technique is to use the Gallai-Edmonds structure
theorem for matchings.
\begin{thm}
\begin{eqnarray*}
rb(n,kK_2)=
\begin{cases}
4, & n=4 \text{ and } k=2;\\
ext(n,(k-1)K_2)+3, & n=2k \text{ and } k\geq7;\\
ext(n,(k-1)K_2)+2, & \text{otherwise.}
\end{cases}
\end{eqnarray*}
\end{thm}

\section{Preliminaries}

Let $M$ be a matching in a given graph $G$. Then the subgraph of $G$
induced by $M$, denoted by $\langle M\rangle_G$ or $\langle
M\rangle$, is the subgraph of $G$ whose edge set is $M$ and whose
vertex set consists of the vertices incident with some edges in $M$.
A vertex of $G$ is said to be \emph{saturated} by $M$ if it is
incident with an edge of $M$; otherwise, it is said to be
\emph{unsaturated}. If every vertex of a vertex subset $U$ of $G$ is
saturated, then we say that $U$ is saturated by $M$. A matching with
maximum cardinality is called a maximum matching.

In a given graph $G$, $N_G(U)$ denotes the set of vertices of $G$
adjacent to a vertex of $U$. If $R,T\in V(G)$, we denote $E_G(R,T)$
or $E(R,T)$ as the set of all edges having a vertex from both $R$
and $T$. Let $G(m,n)$ denote a bipartite graph with bipartition
$A\cup B$, and $|A|=m$ and $|B|=n$. Without loss of generality, in
the following we always assume that $m\geq n$.

Let $ext(m,n,H)$ denote the maximum number of edges that a bipartite
graph $G(m,n)$ can have with no subgraph isomorphic to $H$. The
following lemma is due to Ore and can be found in \cite{lov}.
\begin{lem}\label{ore}
Let $G(m,n)$ be a bipartite graph with bipartition $A\cup B$, and
$M$ a maximum matching in $G$. Then the size of $M$ is $m-d$, where
$$d=\max\{|S|-|N_G(S)|: S\subseteq A\}.$$
\end{lem}

We now determine the value $ext(m,n,H)$ for $H=kK_2$.
\begin{thm}\label{ext}
$$ext(m,n,kK_2)=m(k-1)\text{ for all } n\geq k\geq1,$$
that is, for any given bipartite graph $G(m,n)$, if
$|E(G(m,n))|>m(k-1)$, then $kK_2\subset G(m,n)$.
\end{thm}
\pf Suppose that $G$ contains no $kK_2$. Let $M$ be a maximum
matching of $G$ and the size of $M$ is $k-i$, where $i\geq 1$. By
Lemma \ref{ore}, there exists a subset $S\subset A$ such that
$|S|-|N_G(S)|=m-k+i$. Thus
$$|E(G)|\leq|S||N_G(S)|+n(m-|S|)=(|N_G(S)|+m-k+i)|N_G(S)|+n(k-i-|N_G(S)|).$$
Since $0\leq|N_G(S)|\leq k-i\leq k-1$, we obtain
$$|E(G)|\leq \max\{m(k-1),n(k-1)\}= m(k-1).$$
So, $ext(m,n,kK_2)=m(k-1)$.\qed

\begin{lem}\label{pre}
\begin{eqnarray*}
&&ext(2k,(k-1)K_2)=
\begin{cases}
{k-2 \choose 2}+(k-2)(k+2), & 2\leq k\leq7\text;\\
{2k-3 \choose 2}, & k=2 \text{ or }k\geq7.
\end{cases}\\
\end{eqnarray*}
\end{lem}

\pf From Theorem 1.1, we have that $ext(2k,(k-1)K_2)=\max\{{2k-3
\choose 2}, {k-2 \choose 2}+(k-2)(k+2)\}$. Since ${2k-3 \choose
2}-({k-2 \choose 2}+(k-2)(k+2)) =\frac{1}{2}(k-2)(k-7)$, we have
that if $2\leq k\leq7$, $ext(2k,(k-1)K_2)={k-2 \choose
2}+(k-2)(k+2)$, and if
$k=2$ or $k\geq7$, $ext(2k,(k-1)K_2)={2k-3 \choose 2}$.\qed\\

Let $G$ be a graph. Denote by $D(G)$ the set of all vertices in $G$
which are not covered by at least one maximum matching of $G$. Let
$A(G)$ be the set of vertices in $V(G)-D(G)$ adjacent to at least
one vertex in $D(G)$. Finally let $C(G)=V(G)-A(G)-D(G)$. We denote
the $D(G)$, $A(G)$ and $C(G)$ as the {\it canonical decomposition}
of $G$.

A {\it near-perfect} matching in a graph $G$ is a matching of $G$
covering all but exactly one vertex of $G$. A graph $G$ is said to
be {\it factor-critical} if $G-v$ has a perfect matching for every
$v\in V(G)$.

\begin{thm} (The Gallai-Edmonds Structure Theorem
\cite{lov})\label{str} \ For a graph $G$,  let $D(G)$, $A(G)$ and
$C(G)$ be defined as above. Then
\begin{itemize}
\item[(a)] The components of the subgraph induced by $D(G)$ are
factor-critical.
\item[(b)] The subgraph induced by $C(G)$ has a perfect matching.
\item[(c)] The bipartite graph obtained from $G$ by deleting the
vertices of $C(G)$ and the edges spanned by $A(G)$ and by
contracting each component of $D(G)$ to a single vertex has positive
surplus (as viewed from $A(G)$).
\item[(d)] Any maximum matching $M$ of $G$ contains a
near-perfect matching of each component of $D(G)$, a perfect
matching of each component of $C(G)$ and matches all vertices of
$A(G)$ with vertices in distinct components of $D(G)$.
\item[(e)] The size of a maximum matching $M$ is
$\frac{1}{2}(|V(G)|-c(D(G))+|A(G)|)$, where $c(D(G))$ denotes the
number of components of the graph spanned by $D(G)$.\qed
\end{itemize}
\end{thm}

\section{Main results}

For $k=1$, it is clear that $rb(n,K_2)$=1. Now we determine the
value of $rb(n,2K_2)$ (for $k=2$).

\begin{thm}
$$rb(4,2K_2)=4,$$
and
$$rb(n,2K_2)=2=ext(n,K_2)+2 \text{ for all }n\geq5.$$
\end{thm}
\pf It is obvious that $rb(4,2K_2)\leq 4$. Let
$V(K_4)=\{a_1,a_2,a_3,a_ 4\}$. If $K_4$ is edge-colored with 3
colors such that $c(a_1a_2)=c(a_3a_4)=1$, $c(a_1a_3)=c(a_2a_4)=2$
and $c(a_1a_4)=c(a_2a_3)=3$, then $K_4$ contains no TMC $2K_2$. So,
$rb(4,2K_2)=4$.

For $n\geq5$, let the edges of $G=K_n$ be colored with at least 2
colors. Suppose that $K_n$ contains no TMC $2K_2$. Let $e_1=a_1b_1$
be an edge with $c(e_1)=1$, $T=\{a_1,b_1\}$ and $R=V(K_n)-T$. Then
$c(e)=1$ for all edges $e\in E(G[R])$. Moreover, $c(e)=1$ for all
edges $e\in E(T,R)$, since $|R|\geq3$. But then $K_n$ is
monochromatic, a contradiction. So, $rb(n,2K_2)=2$ for all
$n\geq5$.\qed

The next proposition provides a lower and upper bound for
$rb(n,kK_2)$.

\begin{prop} \
$ext(n,(k-1)K_2)+2\leq rb(n,kK_2)\leq ext(n,kK_2)+1$.
\end{prop}
\pf The upper bound is obvious. For the lower bound, an extremal
coloring of $K_n$ can be obtained from an extremal graph $S_n$ for
$ext(n,(k-1)K_2)$ by coloring the edges of $S_n$ differently and the
edges of $\overline{S_n}$ by
one extra color. It is obvious that the coloring does not contain a TMC $kK_2$.\qed\\

We will show that the lower bound can be achieved for all $n\geq
2k+1$ and $k\geq3$, and thus obtain the exact value of $rb(n,kK_2)$
for all $n\geq2k+1$ and $k\geq3$.

For $n=2k$, we suppose that $H=K_{2k-3}$ is a subgraph of $K_n$ and
$V(K_n)-V(H)=\{a_1,a_2,a_3\}$. If $K_n$ is edge-colored such that
$c(a_1a_2)=1$, $c(a_1a_3)=c(a_2a_3)=2$, $c(e)=1$ for all edges $e\in
E(a_3,V(H))$, $c(e)=2$ for all edges $e\in E(a_1,V(H))\cup
E(a_2,V(H))$ and the edges of $H=K_{2k-3}$ is colored differently by
${2k-3\choose2}$ extra colors. It is easy to check that the coloring
does not contain a TMC $kK_2$ in $K_n$. So, $rb(2k,kK_2)\geq
{2k-3\choose2}+3$ for all $k\geq3$. Hence, if $k\geq7$, then
$ext(2k,(k-1)K_2)={2k-3\choose2}$ and $rb(2k,kK_2)\geq
ext(2k,(k-1)K_2)+3$. We will show that the lower bound can be
achieved for all $n\geq 2k$ and $k\geq7$.

\begin{thm}\label{thm1}
For all $n\geq2k$ and $k\geq3$, we have
\begin{equation*}
rb(n,kK_2)=
\begin{cases}
ext(n,(k-1)K_2)+3,  & n=2k \text{ and } k\geq7\text{;}\\
ext(n,(k-1)K_2)+2,  & \text{ otherwise. }
\end{cases}
\end{equation*}
\end{thm}

\pf We shall prove the theorem by contradiction. If $n=2k$ and
$k\geq7$, let the edges of $K_n$ be colored with $ext(n,(k-1)K_2)+3$
colors; otherwise, let the edges of $K_n$ be colored with
$ext(n,(k-1)K_2)+2$ colors. Suppose that $K_n$ contains no TMC
$kK_2$. Now let $G\subset K_n$ be a TMC spanning subgraph which
contains all colors in $K_n$, i.e., if $n=2k$ and $k\geq7$,
$|E(G)|=ext(n,(k-1)K_2)+3$; otherwise $|E(G)|=ext(n,(k-1)K_2)+2$.
Since $|E(G)|\geq ext(n,(k-1)K_2)+2$, there is a TMC $(k-1)K_2$ in
$G$.

First, we prove the following two assertions.\\

\noindent{\textbf{Claim 1}}: If two components of $G$ consist of a
$K_{2k-3}$ and a $K_3$, respectively, and the other components are
isolated vertices (see Figure 1), then $K_n$ contains a TMC $kK_2$.

Denote $SG_1$ as the special graph $G$ and $Q$ as the set of
isolated vertices of $G$. Without loss of generality, we suppose
that $c(u_1u_2)=1, c(u_2u_3)=2, c(u_1u_3)=3, c(v_1v_2)=4,
c(v_2v_3)=5, c(v_1v_3)=6$ (see Figure 1). The proof of the claim is
given by distinguishing the following two cases:
\begin{figure}
\begin{center}
\psfrag{u1}{$u_1$}\psfrag{u2}{$u_2$}\psfrag{u3}{$u_3$}
\psfrag{v1}{$v_1$}\psfrag{v2}{$v_2$}\psfrag{v3}{$v_3$}
\psfrag{K_{2k-3}}{$K_{2k-3}$}\psfrag{K_3}{$K_3$}
\psfrag{Q}{$Q$}
\includegraphics[width=5cm]{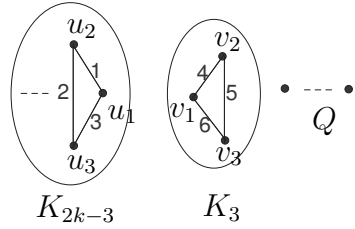}
\caption{The special graph $SG_1$.}
\end{center}
\end{figure}

Case I. $k\geq4$.

We suppose that $G$ contains no TMC $kK_2$. We will show
$c(u_1v_1)=5$. If $c(u_1v_1)\neq5$, then in $G_1=K_{2k-3}-u_1$ the
number of edges whose colors are not $c(u_1v_1)$ is at least
${2k-4\choose2}-1$. Since $k\geq 4$, we have
${2k-4\choose2}-1>ext(2k-4,(k-2)K_2)={2k-5\choose2}$. Thus we can
obtain a TMC $H=(k-2)K_2$ which contains no color $c(u_1v_1)$ in
$G_1$, and hence there is a TMC $kK_2=H\cup \{u_1v_1,v_2v_3\}$ in
$K_n$. So, $c(u_1v_1)$ must be 5. By the same token, $c(u_2v_2)$ and
$c(u_3v_3)$ must be 6 and 4, respectively. Now we can obtain a TMC
$H'=(k-3)K_2$ in $G_2=K_{2k-3}-u_1-u_2-u_3$, and hence there is a
TMC $kK_2=H'\cup \{u_1v_1,u_2v_2,u_3v_3\}$ in $K_n$.

Case II. $k=3$.

We suppose that $K_n$ contains no TMC $3K_2$. Then
$c(u_1v_1)\in\{2,5\}, c(u_2v_2)\in\{3,6\}, c(u_3v_3)\in\{1,4\}$. Now
we can obtain a TMC $3K_2=u_1v_1\cup u_2v_2\cup u_3v_3$ in
$K_n$.\\

\noindent{\textbf{Claim 2}}: If $n\geq2k+1$ and two components of
$G$ are $G'$ and $G''$, where $G'$ and $G''$ is a $K_{2k-3}$ and a
$P_3$, respectively, or $G'$ and $G''$ is a $K^-_{2k-3}$ and a
$K_3$, respectively, and the other components are isolated vertices
(see Figure 2), then $K_n$ contains a TMC $kK_2$, where $P_3$ is a
path with three vertices and $K^-_{2k-3}$ is obtained from
$K_{2k-3}$ by deleting an edge.

Denote $SG_2$ as the special graph $G$ and $Q$ as the set of
isolated vertices of $G$. Without loss of generality, we suppose
that $c(u_1u_2)=1, c(u_2u_3)=2, c(u_1u_3)=3, c(v_1v_2)=4,
c(v_2v_3)=5$ (see Figure 2). The proof of the claim is given by
distinguishing the following two cases:

\begin{figure}
\begin{center}
\psfrag{u1}{$u_1$}\psfrag{u2}{$u_2$}\psfrag{u3}{$u_3$}
\psfrag{v1}{$v_1$}\psfrag{v2}{$v_2$}\psfrag{v3}{$v_3$}
\psfrag{v4}{$v_4$}
\psfrag{K_{2k-3} or K^-_{2k-3}}{$K_{2k-3}$ or
$K^-_{2k-3}$} \psfrag{K_3 or P_3}{$K_3$ or $P_3$} \psfrag{Q}{$Q$}
\psfrag{G'}{$G'$}\psfrag{G''}{$G''$}
\includegraphics[width=6.5cm]{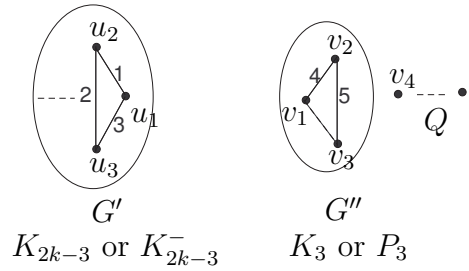}
\caption{The special graph $SG_2$. $G'$ and $G''$ is a $K_{2k-3}$
and a $P_3$, respectively, or $G'$ and $G''$ is a $K^-_{2k-3}$ and a
$K_3$, respectively.}
\end{center}
\end{figure}

Case I. $k\geq4$.

Since $n\geq2k+1$, we suppose that $v_4\in Q$. If $c(u_1v_4)=j$,
without loss of generality, we suppose that $j\neq4$. The number of
edges of $G'-u_1$ whose color is not $j$ is at least ${2k-4\choose
2}-2$ and ${2k-4\choose2}-2>ext(2k-4,(k-2)K_2)={2k-5\choose2}$. Then
there is a TMC $H=(k-2)K_2$ in $G'-u_1$ which contains no color $j$.
We can obtain a TMC $kK_2=H\cup u_1v_4\cup v_1v_2$ in $K_n$.

Case II. $k=3$.

Without loss of generality, we suppose that $G'$ and $G''$ is a
$K_3$ and a $P_3$, respectively. We suppose that $K_n$ contains no
TMC $3K_2$. Then, $c(u_1v_4)\in\{2,5\}\cap\{2,4\}$, i.e.,
$c(u_1v_4)=2$, $c(u_3v_3)\in\{2,4\}\cap\{1,4\}$, i.e.,
$c(u_1v_4)=4$, $c(u_2v_1)\in\{2,5\}\cap\{3,5\}$, i.e.,
$c(u_1v_4)=5$. Now we obtain a TMC
$3K_2=u_1v_4\cup u_3v_3\cup u_2v_1$. See Figure 3.\\

\begin{figure}
\begin{center}
\psfrag{u1}{$u_1$}\psfrag{u2}{$u_2$}\psfrag{u3}{$u_3$}
\psfrag{v1}{$v_1$}\psfrag{v2}{$v_2$}\psfrag{v3}{$v_3$}
\psfrag{v4}{$v_4$}
\psfrag{K_3}{$K_3$} \psfrag{P_3}{$P_3$}
\psfrag{Q}{$Q$}
\includegraphics[width=6cm]{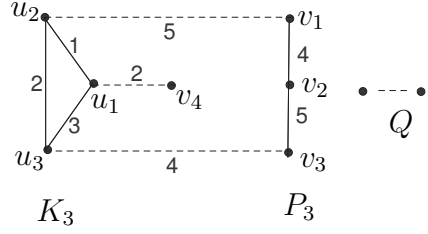}
\caption{We can obtain a TMC $3K_2=u_1v_4\cup u_3v_3\cup u_2v_1$ in
$K_n$.}
\end{center}
\end{figure}

Let $D(G)$, $A(G)$, $C(G)$ as the {\it canonical decomposition} of
$G$ and $c(D(G))=q$, $|A(G)|=s$, $|V(G)|=n$. Since the size of the
maximum matchings of $G$ is $k-1$, by Theorem \ref{str} $(e)$,
$k-1=\frac{1}{2}(n-q+s)$, i.e., $q=n-2k+2+s$. Let the components of
$D(G)$ be $D_1$, $D_2$, $\cdots$, $D_q$. By Theorem \ref{str} $(a)$,
the components of the subgraph induced by $D(G)$ are
factor-critical, hence we suppose that  $|V(D_i)|=2l_i+1$ for $1\leq
i\leq q$, without loss of generality, $l_1\geq l_2\geq\cdots\geq
l_q\geq0$. Let the components of $C(G)$ be $C_1$, $C_2$, $\cdots$,
$C_{q'}$ with $|V(C_i)|=2t_i$ for $1\leq i\leq q'$.

Since $s+q=s+n-2k+2+s\leq n$, then $0\leq s\leq k-1$. Moreover,
\begin{eqnarray*}
n=s+\sum_{i=1}^q(2l_i+1)+|C(G)|&\geq &
s+(2l_1+1)+\sum_{i=2}^q(2l_i+1)\\
&\geq & s+(2l_1+1)+(q-1)\\
&\geq &s+(2l_1+1)+(n-2k+2+s-1),
\end{eqnarray*}
 hence
$2l_1+1\leq2k-2s-1.$

Now we distinguish four cases to finish the proof of the theorem.\\

Case 1. $s=k-1$.

In this case, since $s+q=(k-1)+n-2k+2+(k-1)=n$, then
$C(G)=\emptyset$ and $l_1=l_2=\cdots=l_q=0$. The components of the
subgraph induced by $D(G)$ are isolated vertices. We distinguish two
subcases to finish the proof of the case.

Subcase 1.1. There is at most one vertex $u$ in $D(G)$ such that
$d_G(u)<k-1$.

We suppose $v\in D(G)$ and $u\neq v$. Let $G(n-k-1,k-1)$ be the
bipartite graph obtained from $G$ by deleting the vertices $u,v$ and
the edges spanned by $A(G)$. It is obvious that $uv\in E(K_n)$ and
$uv\notin E(G)$, without loss of generality, we suppose $c(uv)=1$.
Then the number of edges in $G(n-k-1,k-1)$ whose color is not 1 is
at least $(n-k-1)(k-1)-1$. Since $n-k-1\geq2$, then $(n-k-1)(k-1)-1>
ext(n-k-1,k-1,(k-1)K_2)=(n-k-1)(k-2)$. By Lemma \ref{ext}, there
exists a TMC $H=(k-1)K_2$ in $G(n-k-1,k-1)$ which contains no color
1, thus we obtain a TMC $kK_2=H\cup uv$ in $K_n$.

Subcase 1.2. There exist at least two vertices $u$, $v$ in $D(G)$
such that $d_G(u)<k-1$ and $d_G(v)<k-1$.

We suppose that $c(uv)=1$. Let $G'(n-k-1,k-1)$ be the bipartite
graph obtained from $G$ by deleting the vertices $u,v$ and the edges
spanned by $A(G)$ and the edge whose color is 1. Thus there is no
TMC $(k-1)K_2$ in $G'(n-k-1,k-1)$. Hence, by Lemma \ref{ext},
\begin{eqnarray*}
|E(G)|&\leq&1+ext(n-k-1,k-1,(k-1)K_2)+2(k-2)+{k-1\choose2}\\
&\leq&1+(k-2)(n-k-1)+2(k-2)+{k-1\choose2}\\
&=&{k-2\choose2}+(k-2)(n-k+2)+1\\
&<&ext(n,(k-1)K_2)+2,
\end{eqnarray*}
which contradicts $|E(G)|\geq ext(n,(k-1)K_2)+2$.\\

Case 2. $0\leq s\leq k-2$ and $2l_1+1\leq2k-2s-3$.

In this case, if $2k-2s-3=1$, then $l_1=l_2=\cdots=l_q=0$, $s=k-2$
and $|C(G)|=2$, hence
\begin{eqnarray*}
|E(G)|&\leq&{s\choose2}+s(n-s)+{2\choose2}\\
&=&{k-2\choose2}+(k-2)(n-k+2)+1\\
&<&ext(n,(k-1)K_2)+2,
\end{eqnarray*}
which contradicts $|E(G)|\geq ext(n,(k-1)K_2)+2$.

If $2k-2s-3\geq3$, then $0\leq s\leq k-3$ and
\begin{eqnarray*}
\sum_{i=2}^q(2l_i+1)+\sum_{i=1}^{q'}(2t_i)&=&n-s-(2l_1+1)\\
&\geq&n-s-(2k-2s-3)=(q-1)+2.
\end{eqnarray*}
Thus, if $|C(G)|\geq2$, then
\begin{eqnarray*}
|E(G)|&\leq&{s\choose2}+s(n-s)+\sum_{i=1}^q{2l_i+1\choose2}+\sum_{i=1}^{q'}{2t_i\choose2}\\
&\leq&{s\choose2}+s(n-s)+{2l_1+1+\sum_{i=2}^q2l_i\choose 2}+\sum_{i=1}^{q'}{2t_i\choose2}\\
&\leq&{s\choose2}+s(n-s)+{2l_1+1+\sum_{i=2}^q2l_i+(\sum_{i=1}^{q'}2t_i-2)\choose 2}+{2\choose2}\\
&=&{s\choose2}+s(n-s)+{n-s-(q-1)-2\choose2}+{2\choose2}\\
&=&{s\choose2}+s(n-s)+{2k-2s-3\choose2}+{2\choose2}:=f_1(s)\\
\end{eqnarray*}
Hence,
\begin{eqnarray*}
f_1(0)&=&{2k-3\choose2}+1<ext(n,(k-1)K_2)+2,\\
f_1(k-3)&=&{k-2\choose2}+(k-2)(n-k+2)-(n-k)+2\\
&<&{k-2\choose2}+(k-2)(n-k+2)<ext(n,(k-1)K_2)+2.
\end{eqnarray*}
Since $0\leq s\leq k-3$,
$|E(G)|\leq\max\{f_1(0),f_1(k-3)\}<ext(n,(k-1)K_2)+2$, which
contradicts $|E(G)|\geq ext(n,(k-1)K_2)+2$.

If $|C(G)|=0$, then $2l_2+1\geq3$ and
\begin{eqnarray*}
|E(G)|&\leq&{s\choose2}+s(n-s)+\sum_{i=1}^q{2l_i+1\choose2}+\sum_{i=1}^{q'}{2t_i\choose2}\\
&\leq&{s\choose2}+s(n-s)+{2l_1+1+\sum_{i=3}^q2l_i+\sum_{i=1}^{q'}2t_i\choose 2}+{2l_2+1\choose2}\\
&\leq&{s\choose2}+s(n-s)+{2l_1+1+\sum_{i=3}^q2l_i+\sum_{i=1}^{q'}2t_i+(2l_2-2)\choose 2}+{3\choose2}\\
&=&{s\choose2}+s(n-s)+{n-s-(q-1)-2\choose2}+{3\choose2}\\
&=&{s\choose2}+s(n-s)+{2k-2s-3\choose2}+{3\choose2}:=f_2(s).\\
\end{eqnarray*}
Thus,
\begin{eqnarray*}
f_2(0)&=&{2k-3\choose2}+3,\\
f_2(1)&=&{2k-3\choose2}+n-4k+11,\\
f_2(k-3)&=&{k-2\choose2}+(k-2)(n-k+2)-(n-k)+4\\
&\leq&{k-2\choose2}+(k-2)(n-k+2)+1<ext(n,(k-1)K_2)+2.
\end{eqnarray*}

If $s=0$ and $|E(G)|={2k-3\choose2}+3$, then $G\cong SG_1$. By Claim
1, we can obtain a TMC $kK_2$ in $K_n$. If $s=0$, $n\geq2k+1$ and
$|E(G)|={2k-3\choose2}+2$, then $G\cong SG_2$. By Claim 2, we can
obtain a TMC $kK_2$ in $K_n$. So, if $n\geq2k+1$, then
$|E(G)|\leq{2k-3\choose2}+1<ext(n,(k-1)K_2)+2$, which contradicts
$|E(G)|=ext(n,(k-1)K_2)+2$. If $n=2k$ and $k\geq7$, then
$|E(G)|\leq{2k-3\choose2}+2=ext(n,(k-1)K_2)+2$, which contradicts
$|E(G)|=ext(n,(k-1)K_2)+3$. If $n=2k$ and $3\leq k\leq6$, then
$|E(G)|\leq{2k-3\choose2}+2\leq
{k-2\choose2}+(k-2)(k+2)=ext(n,(k-1)K_2)$, which contradicts
$|E(G)|=ext(n,(k-1)K_2)+2$.

If $1\leq s\leq k-3$, then $k\geq4$ and
$|E(G)|\leq\max\{f_2(1),f_2(k-3)\}$. So, if $f_2(k-3)\geq f_2(1)$,
then $|E(G)|\leq f_2(k-3)<ext(n,(k-1)K_2)+2$, a contradiction. If
$f_2(1)>f_2(k-3)$, then
$$
{2k-3\choose2}+n-4k+11>{k-2\choose2}+(k-2)(n-k+2)-(n-k)+4.
$$
Hence $2k\leq n<\frac{1}{2}(5k-7)$, $k>7$ and
\begin{eqnarray*}
|E(G)|&\leq&f_2(1)={2k-3\choose2}+n-4k+11\\
&<&{2k-3\choose2}+\frac{1}{2}(15-3k)\\
&<&ext(n,(k-1)K_2)+2,
\end{eqnarray*}
a contradiction.\\

Case 3. $0\leq s\leq k-2$, $2l_1+1=2k-2s-1$ and $n\geq2k+1$.

In this case, $s+(2l_1+1)+(q-1)=n$, hence $C(G)=\emptyset$,
$l_2=l_3=\cdots=l_q=0$ and each $D_i$ for $2\leq i\leq q$ is an
isolated vertex.

Let $G(q,s)$ be the bipartite graph obtained from $G$ by deleting
the edges spanned by $A(G)$ and by contracting the component $D_1$
to a single vertex $p$. Thus by Theorem \ref{str} (c) and (d), we
can obtain a maximum matching $M$ of size $k-1$ such that $M$
contains a maximum matching $M_1$ of $G(q,s)$ which does not match
vertex $p$ and a near-perfect matching $M_2$ of $D_1$. Since
$q=n-2k+2+s\geq s+3$, there exist two vertices $u,v\in D(G)-D_1$ and
$u,v\notin \langle M\rangle$. It is obvious that $uv\in E(K_n)$ and
$uv\notin E(G)$. We suppose that $c(uv)=1$, hence there exists an
edge $e=yz\in M$ with $c(e)=1$. Now we distinguish two subcases to
complete the proof of the case.

Subcase 3.1. $e\in M_1$.

In this subcase, $s\geq1$ and $yz\in E_G(A(G),D(G))$, without loss
of generality, we suppose that $y\in A(G)$. If there exists an edge
$yz_1\in E_G(y,D_1)$ with $z_1\in D_1$, then we can obtain another
maximum matching $M'_1$ of $G(q,s)$ with $M'_1=M_1\cup yz_1-yz$ and
a near-perfect matching $M'_2$ of $D_1$ which does not match $z_1$.
Thus we obtain a TMC $kK_2=M'_1\cup M'_2\cup uv$ in $K_n$. See
Figure 4.
\begin{figure}
\begin{center}
\psfrag{u}{$u$}\psfrag{v}{$v$}\psfrag{y}{$y$}\psfrag{z}{$z$}
\psfrag{p}{$p$}\psfrag{e}{$e$}
\psfrag{z_1}{$z_1$}\psfrag{M_1}{$M_1$}
\psfrag{A(G)}{$A(G)$}\psfrag{s>0}{$s\geq1$} \psfrag{q>s+2}{$q\geq
s+3$}\psfrag{|D_1|=2k-2s-1}{$|D_1|=2k-2s-1$}
\includegraphics[width=7cm]{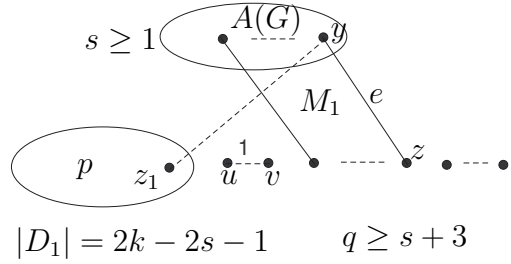}
\caption{If $yz_1\in E_G(y,D_1)$, we can obtain a TMC $kK_2=M'_1\cup
M'_2\cup uv$ in $K_n$.}
\end{center}
\end{figure}

Thus we suppose that $E_G(y,D_1)=\emptyset$. There is no matching of
size $s$ in $G'(q-3,s)=G(q,s)-p-u-v-e$. By Lemma \ref{ext},
$|E_G(G')|\leq(s-1)(q-3)=(s-1)(n-2k+s-1)$. Now
\begin{eqnarray*}
|E(G)|&\leq&{s\choose2}+{2k-2s-1\choose2}+1+|E_G(G')|\\
&&+|E_G(D_1,A(G))|+|E_G(\{u,v\},A(G))|\\
&\leq&{s\choose2}+{2k-2s-1\choose2}+1+(s-1)(n-2k+s-1)\\
&&+(2k-2s-1)(s-1)+2s:=f_3(s)\\
\end{eqnarray*}
Hence,
\begin{eqnarray*}
f_3(1)&=&{2k-3\choose2}+3,\\
f_3(2)&=&{2k-3\choose2}+n-4k+11,\\
f_3(k-2)&=&{k-2\choose2}+(k-2)(n-k+2)-(n-k)+4\\
&\leq&{k-2\choose2}+(k-2)(n-k+2)<ext(n,(k-1)K_2)+2.
\end{eqnarray*}

If $s=1$, then $|E(G)|\leq{2k-3\choose2}+3$. If
$|E(G)|={2k-3\choose2}+3$, then $(G-e+uv)\cong SG_1$. By the proof
of Claim 1, we can obtain a TMC $kK_2$ in $K_n$. If
$|E(G)|={2k-3\choose2}+2$, then $(G-e+uv)\cong SG_2$. By the proof
of Claim 2 , we can obtain a TMC $kK_2$ in $K_n$. If
$|E(G)|\leq{2k-3\choose2}+1\leq ext(n,(k-1)K_2)+1$, this contradicts
$|E(G)|=ext(n,(k-1)K_2)+2$.

If $2\leq s\leq k-2$, then $k\geq4$ and
$|E(G)|\leq\max\{f_3(2),f_3(k-2)\}$. So, if $f_3(k-2)\geq f_3(2)$,
then $|E(G)|\leq f_3(k-2)<ext(n,(k-1)K_2)+2$, a contradiction. If
$f_3(1)>f_3(k-3)$, then
$$
{2k-3\choose2}+n-4k+11>{k-2\choose2}+(k-2)(n-k+2)-(n-k)+4.
$$
Hence, $2k\leq n<\frac{1}{2}(5k-7)$, $k>7$ and
\begin{eqnarray*}
|E(G)|&\leq&f_3(2)={2k-3\choose2}+n-4k+11\\
&<&{2k-3\choose2}+\frac{1}{2}(15-3k)\\
&<&ext(n,(k-1)K_2)+2,
\end{eqnarray*}
a contradiction.

Subcase 3.2. $e\in M_2$.

In this subcase, $y\in D_1$ and $z\in D_1$. By Theorem \ref{str}
(a), $D_1$ is factor-critical, there exists a near-perfect matching
$M'_2$ which does not match $y$, So $M'_2$ does not contain $e=yz$.
Now we obtain a TMC $kK_2=M'_2\cup M_1\cup uv$ in $K_n$.\\

Case 4. $0\leq s\leq k-2$, $2l_1+1=2k-2s-1$ and $n=2k$.

In this case, $q=s+2$ and $s+(2l_1+1)+(q-1)=2k$, hence
$C(G)=\emptyset$, $l_2=l_3=\cdots=l_q=0$ and each $D_i$ for $2\leq
i\leq q$ is an isolated vertex. Now we distinguish two subcases to
complete the proof of the case.

Subcase 4.1. $1\leq s\leq k-2$.

If $E_G(D_1,A(G))=\emptyset$, then
$$|E(G)|\leq{2k-2s-1\choose2}+{s\choose2}+s(s+1):=f_4(s).$$
Thus,
\begin{eqnarray*}
f_4(1)&=&{2k-3\choose2}+2,\\
f_4(k-2)&=&{k-2\choose2}+(k-2)(k+2)+3-3(k-2)\\
\end{eqnarray*}

Since $k\geq3$, then $f_4(1)\geq f_4(k-2)$ and
$|E(G)|\leq\max\{f_4(1),f_4(k-2)\}=f_4(1)={2k-3\choose2}+2$. If
$k\geq7$, this contradicts
$|E(G)|=ext(2k,(k-1)K_2)+3={2k-3\choose2}+3$. If $3\leq k\leq6$,
then
\begin{eqnarray*}
|E(G)|&\leq&{2k-3\choose2}+2\\
&\leq&{k-2\choose2}+(k-2)(k+2)=ext(2k,(k-1)K_2),
\end{eqnarray*}
which contradicts $|E(G)|=ext(2k,(k-1)K_2)+2$.

So we suppose that $E_G(D_1,A(G))\neq\emptyset$. Let $G(s+2,s)$ be
the bipartite graph obtained from $G$ by deleting the edges spanned
by $A(G)$ and by contracting the component $D_1$ to a single vertex
$p$. Thus by Theorem \ref{str} (d), we can obtain a maximum matching
$M$ of size $k-1$ such that $M$ contains a near-perfect matching
$M_1$ of $D_1$ which does not match $w$ with $w\in D_1$ and a
matching $M_2$ of size $s$ which matches all vertices of $A(G)$ with
vertices in $\{w\}\cup (D(G)-D_1)$. Since
$E_G(D_1,A(G))\neq\emptyset$, we can suppose that $w\in \langle
M_2\rangle$. There exist exactly two vertices $u,v\in D(G)-D_1$ and
$u,v\notin \langle M\rangle$. It is obvious that $uv\in E(K_n)$ and
$uv\notin E(G)$. We suppose that $c(uv)=1$, hence there exists an
edge $e=yz\in M$ with $c(e)=1$. Now we distinguish two subcases to
complete the proof of the subcase 4a.

Subcase 4.1.1. $e=yz\in M_1$.

If $s=1$, then $|D_1|=2k-3$ and we suppose $A(G)=\{x\}$. Thus the
size of $M_1$ is $k-2$ and there is no $H=(k-2)K_2$ in
$D'_1=D_1-w-yz$, for otherwise, we can obtain a TMC $kK_2=H\cup
xw\cup uv$ in $K_{2k}$. If $E_G(x,\{y,z\})\neq\emptyset$, say $xy\in
E(G)$, then we can obtain a perfect matching $M'_1$ of $D_1-y$ and a
TMC $kK_2=M'_1\cup uv\cup xy$ in $K_{2k}$. So,
$E_G(x,\{y,z\})=\emptyset$ and
\begin{eqnarray*}
|E(G)|&=&1+|E_G(D'_1)|+|E_G(w,D'_1)|+|E_G(x,D_1)|+|E_G(x,\{u,v\})|\\
&\leq&1+ext(2k-4,(k-2)K_2)+(2k-4)+(2k-5)+2\\
&=&{2k-5\choose2}+4k-6\\
&=&{2k-3\choose2}+3.
\end{eqnarray*}

Denote $SG_3$ be the special graph $G$ shown in Figure 5, whence
$E(SG_3)=E(K^-_{2k-3})\cup xu\cup xv\cup yw\cup yz$.  Without loss
of generality, we suppose that $c(wy)=4$. If
$|E(G)|={2k-3\choose2}+3$, it is easy to check that $G\cong SG_3$.

If $k\geq7$, then by the beginning hypothesis
$|E(G)|=ext(2k,(k-1)K_2)+3={2k-3\choose2}+3$, whence $G\cong SG_3$.
Now ${2k-4\choose2}-1>ext(2k-4,(k-2)K_2)$, we can obtain a TMC
$H=(k-2)K_2$ in $K^-_{2k-3}-w$, whence a TMC $kK_2=H\cup yw\cup uv$
in $K_{2k}$. If $3\leq k\leq6$, then
$${2k-3\choose2}+3\leq{k-2\choose2}+(k-2)(k+2)+1=ext(2k,(k-1)K_2)+1,$$
which contradicts $|E(G)|=ext(2k,(k-1)K_2)+2$.

\begin{figure}
\begin{center}
\psfrag{u}{$u$}\psfrag{v}{$v$}\psfrag{w}{$w$}
\psfrag{x}{$x$}\psfrag{y}{$y$}\psfrag{z}{$z$}
\psfrag{K_{2k-3}-xz}{$k^-_{2k-3}=K_{2k-3}-xz$} \psfrag{SG_3}{$SG_3$}
\includegraphics[width=5cm]{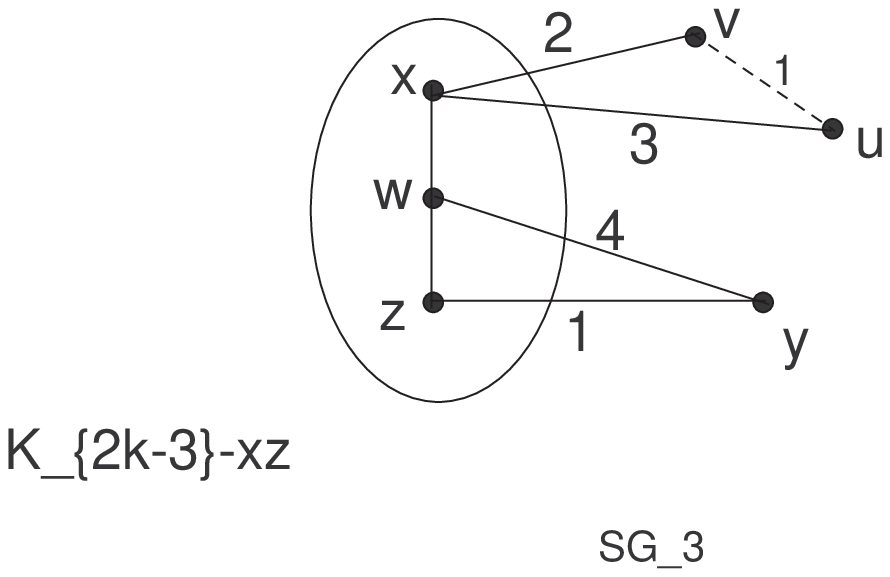}
\caption{The special graph $SG_3$ and $|E(SG_3)|={2k-3\choose2}+3$.}
\end{center}
\end{figure}

If $2\leq s\leq k-2$, then $k\geq4$. We suppose that $x\in A(G)$ and
$xw\in M_2$. By the the same token, $E_G(x,\{y,z\})=\emptyset$ and
there is no $(k-s-1)K_2$ in $D'_1=D_1-w-yz$. If
$E_G(A(G)-x,\{y,z\})\neq\emptyset$, say $x'y\in E(G)$, then there is
no $H=(s-1)K_2$ in bipartite graph $G'(s-1,s-1)=G-\{D_1\cup u\cup
v\cup x'\}$, for otherwise, we can obtain a perfect matching $M'_1$
in $D_1-y$ and a TMC $kK_2=M'_1\cup H\cup uv\cup x'y$. See Figure 6.
Thus,
\begin{eqnarray*}
|E_G(A(G),D(G))|&=&|E_G(A(G),D_1-y-z)|+|E(A(G),\{y,z\})|\\
&&+|E_G(A(G),\{u,v\})|+|E_G(G'(s-1,s-1))|\\
&&+|E_G(x',D(G)-D_1-u-v)|\\
&\leq&(2k-2s-3)s+2(s-1)+2s\\
&&+ext(s-1,s-1,(s-1)K_2)+(s-1)\\
&=&(2k-2s-3)s+2s+(s-1)(s+1).
\end{eqnarray*}

\begin{figure}
\begin{center}
\psfrag{u}{$u$}\psfrag{v}{$v$}\psfrag{w}{$w$}
\psfrag{x}{$x$}\psfrag{y}{$y$}\psfrag{z}{$z$} \psfrag{e}{$e$}
\psfrag{x'}{$x'$}\psfrag{D_1}{$D_1$}\psfrag{M_2}{$M_2$}
\psfrag{A(G)}{$A(G)$} \psfrag{q=s+2}{$q=s+2$}\psfrag{s-1}{$s-1$}
\includegraphics[width=6cm]{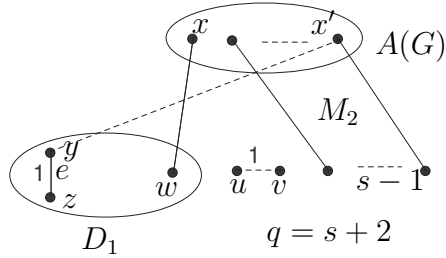}
\caption{There is no $(k-s-1)K_2$ in $D'_1=D_1-w-yz$. If $x'y\in
E(G)$, there is no $(s-1)K_2$ in bipartite graph
$G'(s-1,s-1)=G-\{D_1\cup u\cup v\cup x'\}$.}
\end{center}
\end{figure}
If $E_G(A(G)-x,\{y,z\})=\emptyset$, then
\begin{eqnarray*}
|E_G(A(G),D(G))|&=&|E_G(A(G),D_1-y-z)|+|E_G(A(G),D(G)-D_1)|\\
&\leq&(2k-2s-3)s+s(s+1).
\end{eqnarray*}
So,
\begin{eqnarray*}
&&|E_G(A(G),D(G))|\\
&\leq&\max\{(2k-2s-3)s+2s+(s-1)(s+1),(2k-2s-3)s+s(s+1)\}\\
&=&(2k-2s-3)s+2s+(s-1)(s+1).
\end{eqnarray*}
Now, we have
\begin{eqnarray*}
|E(G)|&=&{s\choose2}+1+|E_G(D'_1)|+|E_G(w,D'_1)|+|E_G(A(G),D(G))|\\
&\leq&{s\choose2}+1+{2k-2s-3\choose2}+(2k-2s-2)+\\
&&(2k-2s-3)s+2s+(s-1)(s+1):=f_5(s).
\end{eqnarray*}
Thus,
\begin{eqnarray*}
f_5(2)&=&{2k-3\choose2}-2k+11,\\
f_5(k-2)&=&{k-2\choose2}+(k-2)(k+2)-k+4\\
&<&ext(2k,(k-1)K_2)+2.
\end{eqnarray*}

If $4\leq k\leq6$, then $f_5(k-2)\geq f_5(2)$ and
$|E(G)|\leq\max\{f_5(2),f_5(k-2)\}=f_5(k-2)<ext(2k,(k-1)K_2)+2$,
which contradicts $|E(G)|=ext(2k,(k-1)K_2)+2$.

If $k\geq7$, then $f_5(2)\geq f_5(k-2)$ and
$|E(G)|\leq\max\{f_5(2),f_5(k-2)\}=f_5(2)={2k-3\choose2}-2k+11<{2k-3\choose2}=ext(2k,(k-1)K_2),$
which contradicts $|E(G)|=ext(2k,(k-1)K_2)+3$.

Subcase 4.1.2. $e=yz\in M_2$.

Without loss of generality, we suppose that $y\in A(G)$.

If $s=1$, then $A(G)=\{y\}$, $yz=yw$ and $c(yw)=c(uv)=1$. Then
$E_G(y,D_1-w)=\emptyset$, for otherwise, say $yw'\in E_G(y,D_1-w)$
with $w'\in (D_1-w)$, we can obtain a TMC $H=(k-2)K_2$ in $D_1-w'$
and a TMC $kK_2=H\cup yw'\cup uv$ in $K_{2k}$. So,
\begin{eqnarray*}
|E(G)|&=&|E_G(D_1)|+|E_G(y,\{w,u,v\})|\leq{2k-3\choose2}+3.
\end{eqnarray*}
If $3\leq k\leq6$, then
$${2k-3\choose2}+3\leq{k-2\choose2}+(k-2)(k+2)+1=ext(2k,(k-1)K_2)+1,$$
which contradicts $|E(G)|=ext(2k,(k-1)K_2)+2$.

If $k\geq7$, since $|E(G)|={2k-3\choose2}+3$, it is easy to check
that $(G-e+uv)\cong SG_1$. By the proof of the Claim 1, we can
obtain a TMC $kK_2$ in $K_{2k}$.

If $2\leq s\leq k-2$, first we look at the bipartite graph
$G(s+2,s)$. We suppose that $M'_2$ is any maximum matching of size
$s$ in $G(s+2,s)$ with $p\in \langle M'_2\rangle$ and
$u_1,v_1\notin\langle M'_2\rangle$. By Subcase 4.1.1, we can suppose
that there exists an edge $e_1\in M'_2$ such that
$c(e_1)=c(u_1v_1)$. If $d_{G(s+2,s)}(p)=s$ and there is at most one
vertex $u_2$ in $D(G)-D_1$ such that $d_{G(s+2,s)}(u)\leq s-1$, we
suppose $v_2\in D(G)-D_1$ and $u_2\neq v_2$. Let $G(s,s)$ be the
bipartite graph obtained from $G(s+2,s)$ by deleting the vertices
$u_2,v_2$. It is obvious that $u_2v_2\in E(K_n)$ and $u_2v_2\notin
E(G)$. Then the number of edges in $G(s,s)$ whose color is not
$c(u_2v_2)$ is at least $s^2-1$. Since $s\geq2$, then $s^2-1\geq
ext(s,s,sK_2)=s(s-1)+1$. By Lemma \ref{ext}, there exists a TMC
$H=sK_2$ in $G(s,s)$ which contains no color $c(u_2v_2)$, thus we
obtain a TMC $(s+1)K_2=H\cup u_2v_2$. By Theorem \ref{str}, we can
obtain a TMC $kK_2$ in $K_{2k}$.

So, if $d_{G(s+2,s)}(p)=s$, then we suppose there exist at least two
vertices $u_3$, $v_3$ in $D(G)-D_1$ such that $d_{G(s+2,s)}(u_3)\leq
s-1$ and $d_{G(s+2,s)}(v_3)\leq s-1$. Let $G'(s,s)$ be the bipartite
graph obtained from $G(s+2,s)$ by deleting the vertices $u_3,v_3$
and the edge whose color is $c(u_3v_3)$. Thus there is no TMC $sK_2$
in $G'(s,s)$. By Lemma \ref{ext}, $E(G(s+2,s))\leq1+2(s-1)+s(s-1)$
and
$$|E_G(A(G),D(G))|\leq1+2(s-1)+s((2k-2s-1)+(s-2))=1+2(s-1)+s(2k-s-3).$$

Now we suppose that $d_{G(s+2,s)}(p)\leq s-1$. Since
$E(A(G),D_1)\neq\emptyset$, if there exists an edge $w''x'\in
E(A(G),D_1)$ with $x'\in A(G)$, $w''\in D_1$ and $w''x'\neq wx$.
Thus there is no TMC $H=(s-1)K_2$ in $G(s+2,s)-\{p\cup u\cup v\cup
x'\}-yz$, for otherwise, we can obtain a TMC $(s+1)K_2=H\cup uv\cup
w''x'$, a TMC $(k-s-1)K_2$ in $D_1-w''$ and a TMC $kK_2$ in
$K_{2k}$. We have
\begin{eqnarray*}
|E_G(A(G),D(G))|&\leq&|E_G(A(G),D_1)|+(s-1)(s-2)+1\\
&&+|E_G(x',D(G)-D_1-u-v)|+|E_G(A(G),\{u,v\})|\\
&\leq&(2k-2s-1)(s-1)+(s-1)(s-2)+1+(s-1)+2s\\
&=&(2k-2s-1)(s-1)+s^2+2.
\end{eqnarray*}
If $E(A(G),D_1)=\{xw\}$, then
\begin{eqnarray*}
|E_G(A(G),D(G))|&\leq&1+s(s+1).
\end{eqnarray*}

Thus,
\begin{eqnarray*}
&&|E_G(A(G),D(G))|\\
&\leq&\max\{1+2(s-1)+s(2k-s-3),(2k-2s-1)(s-1)+s^2+2,1+s(s+1)\}\\
&=&1+2(s-1)+s(2k-s-3).
\end{eqnarray*}
So,
\begin{eqnarray*}
|E(G)|\leq{s\choose2}+{2k-2s-1\choose2}+1+2(s-1)+s(2k-s-3):=f_6(s).
\end{eqnarray*}
We have
\begin{eqnarray*}
f_6(2)&=&{2k-3\choose2}+3,\\
f_6(3)&=&{2k-3\choose2}-2k+12,\\
f_6(k-2)&=&{k-2\choose2}+(k-2)(k+2)-k+4\\
&<&ext(2k,(k-1)K_2)+2.
\end{eqnarray*}

If $s=2$ and $|E(G)|=f_6(2)={2k-3\choose2}+3$, then it is easy to
check that $G$ has a structure shown in Figure 7. By the proof the
Claim 1, we can obtain a TMC $kK_2$ in $K_{2k}$.
\begin{figure}
\begin{center}
\psfrag{u}{$u$}\psfrag{v}{$v$}\psfrag{K_{2k-3}}{$K_{2k-3}$}
\includegraphics[width=10cm]{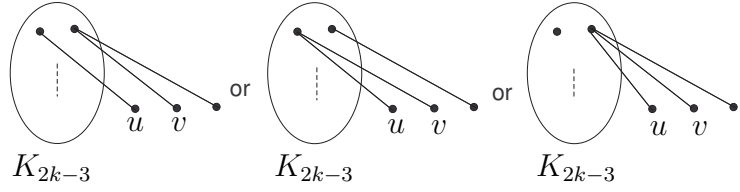}
\caption{$G$ is isomorphic to one of the above three graphs.}
\end{center}
\end{figure}

If $3\leq s\leq k-2$, then $k\geq5$. If $5\leq k\leq6$, then
$f_6(k-2)=f_6(3)$ and $|E(G)|\leq f_6(k-2)<ext(2k,(k-1)K_2)+2$,
which contradicts $|E(G)|=ext(2k,(k-1)K_2)+2$. If $k\geq7$, then
$f_6(3)>f_6(k-2)$ and $|E(G)|\leq
f_6(3)={2k-3\choose2}-2k+12<{2k-3\choose2}=ext(2k,(k-1)K_2)$, which
contradicts $|E(G)|=ext(2k,(k-1)K_2)+3$.\\

Subcase 4.2. $s=0$.

In this subcase, $|V(D_1)|=2k-1$ and $q=2$. We suppose that $z_1\in
D_1$ and $D_2=\{z_2\}$. Let $M$ be a perfect matching of $D_1-z_1$.
Then there exists an edge $e\in M$ such that $c(e)=c(z_1z_2)$. So,
there is no TMC $(k-1)K_2$ in $D_1-z_1-e$. Let $D'_1$ be $D_1-z_1-e$
and $D(D'_1)$, $A(D'_1)$ and $C(D'_1)$ as the {\it canonical
decomposition} of $D'_1$. We look at the graph $G_1=G-e+z_1z_2$. Let
$A'(G_1)=A(D'_1)\cup z_1$ and $D'(G_1)=D(D'_1)\cup z_2$ and
$C'(G_1)=C(D'_1)$. Let $|A'(G_1)|=s'$,
$q'=c(D'(G_1))=c(D(D'_1))+1=(2k-2)-2(k-2)+s-1+1=s+2$. Obviously,
$1\leq s'\leq k-1$. Employing similar technique as in the proofs of
Cases 1, 2 and Subcase 4.1, we can obtain contradictions. The
details are omitted. Up to now, the proof is complete.\qed

\end{document}